\documentstyle[12pt]{article}

\newfont{\frak}{eufm10 scaled\magstep1}
\newfont{\extra}{msbm10 scaled\magstep1}
\newcommand{\fra}[1]{\mbox{\frak #1}}
\newcommand{\extr}[1]{\mbox{\extra #1}}

\newcommand{\sect}[1]{\setcounter{equation}{0}\section{#1}}
\newcommand{\subsect}[1]{\subsection{#1}}
\newcommand{\subsubsect}[1]{\subsubsection{#1}}

\def\be{\begin{equation}}
\def\ee{\end{equation}}
\def\ba{\begin{array}}
\def\ea{\end{array}}
\def\bea{\begin{eqnarray}}
\def\eea{\end{eqnarray}}
\def\C{{\extr C}}
\def\R{{\extr R}}
\def\v{\epsilon}
\def\l{\lambda}
\def\k{\kappa}

\catcode`\@=11
\font\tenmsa=msam10
\font\sevenmsa=msam7
\font\fivemsa=msam5
\font\tenmsb=msbm10
\font\sevenmsb=msbm7
\font\fivemsb=msbm5
\newfam\msafam
\newfam\msbfam
\textfont\msafam=\tenmsa  \scriptfont\msafam=\sevenmsa
  \scriptscriptfont\msafam=\fivemsa
\textfont\msbfam=\tenmsb  \scriptfont\msbfam=\sevenmsb
  \scriptscriptfont\msbfam=\fivemsb

\def\hexnumber@#1{\ifnum#1<10 \number#1\else
 \ifnum#1=10 A\else\ifnum#1=11 B\else\ifnum#1=12 C\else
 \ifnum#1=13 D\else\ifnum#1=14 E\else\ifnum#1=15 F\fi\fi\fi\fi\fi\fi\fi}

\def\msa@{\hexnumber@\msafam}
\def\msb@{\hexnumber@\msbfam}
\mathchardef\blacktriangleright="3\msa@49
\mathchardef\blacktriangleleft="3\msa@4A
\catcode`\@=12

\def\bicross{\triangleright\!\!\!\blacktriangleleft}

\def\leco{\blacktriangleleft}
\def\reco{\blacktriangleright}
\def\act{\triangleright}
\def\RL{\triangleright\!\!\!\blacktriangleleft}
\def\LR{\blacktriangleright\!\!\!\triangleleft}
\def\LLL{\triangleleft}

\def\RIMO{\triangleright\!\!\!<}
\def\LECO{>\!\!\blacktriangleleft}

\begin{document}

\begin{center}
{\large{\bf{Induced Representations of Quantum Groups}}}\footnote
{Talk presented by M.A.del Olmo at the Symposium on Physics and 
Geometry. Zaragoza, February 1998.}
\end{center}

\vskip 1cm

 \begin{center} 
Oscar Arratia $^{1}$ and Mariano A. del Olmo $^{2}$
\end{center}

\begin{center}
{\it $^{1}$ Departamento de Matem\'{a}tica Aplicada a la
Ingenier\'{\i}a,  \\
Universidad de Valladolid. E-47011, Valladolid, Spain.\\
E. mail: oscarr$@$wmatem.eis.uva.es}
\end{center}

\begin{center}
{\it $^{2}$ Departamento de F\'{\i}sica Te\'{o}rica,\\
Universidad de Valladolid. E-47011, Valladolid, Spain.\\ 
E. mail: olmo$@$cpd.uva.es}
\end{center}

\vskip 1cm
\begin{center}
\today
\end{center}

\begin{abstract}
In this paper we show how to construct explicitly induced 
representations for bicrossproduct Hopf algebras with abelian kernels
starting from one-dimensional characters of the commutative sector.
We introduce this technique by means of two concrete physical examples:
two quantum deformations of the $(1+1)$ Galilei algebra. 
\end{abstract}

\sect{Introduction}

From its beginning in 1986 \cite{drinfeld} quantum groups and 
quantum algebras have largely attracted the attention of mathematicians
and physicists. The main reason for this fascination is the very rich
mathematical structure carried by these objects, which allows to mimic 
systematically useful constructions developed in many
other well known  branches of Mathematics, in particular, in Lie group
theory. In consequence, there are a huge variety of potential
applications of quantum groups ranging from integrable systems or
quantum mechanics  to conformal field theory (see, for instance,
Ref. \cite{jimboed} and \cite{gomez}).

The contribution of this paper is located within the applications of
these new mathematical entities to the description of deformed 
symmetries (or $q$--symmetries) of physical systems as well as of the
space-time. Quantum kinematical algebras and groups can
be used for the study of $q$--symmetries of the  deformed space-time 
($q$--space-time), since the $q$--space-time can be considered as a
non-commutative homogeneous space of the quantum kinematical groups.
Consequently, we are interested in quantum groups as the adequate
tool to describe the short range of the space-time structure, which
looks to be non-commutative. 

Another approach to quantum groups is closely related with
the deformation of the commutative algebra of functions.
One of the most interesting examples of it is related with the problem
of the quantization of physical systems and with the
deformation of phase spaces \cite{bayen} (see Ref. \cite{arratia1} for
a review and references therein).   

On the other hand, the study of the representations of the
quantum kinematical groups is an interesting problem, that can be useful
for determining the behaviour of physical systems endowed with deformed
symmetries. Obviously, it looks natural to construct the
representations of quantum groups in the framework of
non-commutative homogeneous spaces, which are the natural arena  for
quantum groups. This procedure fits in Connes' program of
noncommutative geometry \cite{connes}. 

Moreover, it is expected a rich interplay between $q$--spaces and
representations, in particular, in relation with the $q$-analogous of the
harmonic analysis and $q$--special functions. 
  
Physically, as it is well known, a projective unitary irreducible
representation of a symmetry group of a given physical system leads to 
a definition of quantum elementary physical system \cite{wigner}, and
also gives a prescription for computing expected values (the
observables are assumed to form the symmetry algebra we start with).

Kinematical groups like Poincar\'e and Galilei are
semidirect product of the translation group and the homogeneous group of
rotations and boosts (Lorentz or homogeneous Galilei, respectively).
Therefore, the most appropriated method to construct their
unitary representations is the Mackey method for induced representations
of semidirect products \cite{mackey}. 

In this paper we obtain the induced representations of
two non-equivalent quantum deformations of the $(1+1)$ Galilei algebra 
by using a generalization of Mackey's method. In both cases included
here the quantum $(1+1)$ Galilei algebra has a structure of
bicrossproduct, which is  a generalization of the semidirect product of
Lie groups (or algebras) to Hopf algebras \cite{majid}. That
constitutes the first approximation in order to get a quantum analogue
of Mackey's theory. 
Some attempts  have been made to extend this technique to the
quantum case from the mathematical \cite{ibort} as well as from the
physical \cite{giller,maslanka,italianos} point of view. However, in
all these cases the approach has been mainly focused on
corepresentations of quantum groups, in other words, in representations
of the coalgebra part. However, this paper deals with the dual case,
closer to the classical one, constructing representations in the
algebra part. 

The organization of the paper is as follows. In Section 2 we review 
the algebraic structure related with the topics of Hopf algebra,
quantum algebra and quantum group. The bicrossproduct 
structure is also described here. In next Section we introduce the 
basic elements of the theory of induced representations of quantum
groups, which is connected with module theory, and
build up induced representations for two non-equivalent deformations
of the Galilei algebra. Some comments and remarks on the results 
obtained here together with a collection of open problems close the
paper.

\sect{Quantum groups and quantum algebras}

As it is well known quantum groups and quantum algebras are neither
Lie groups nor Lie algebras, but the mathematical structure underlying
both kind of objects is that of Hopf algebra. 

\subsect{Hopf algebras}

A Hopf algebra restores, in some sense, the symmetry lost
when a product law is added to a (complex) vector space $V$ in order to
get an algebra. The Hopf algebraic
setting allows not only for the possibility to compose but also to 
``decompose'' elements in
$V$. More explicitly, on the linear space  $V$ we have two linear
mappings 
$$ 
m:V\otimes V \longrightarrow V, \qquad 
\Delta: V \longrightarrow V \otimes V,
$$
referred to as the product and the coproduct, respectively. Both 
mappings are compatible in the sense that
$$
\Delta \circ m = (m\otimes m) \circ({\rm id}\circ \tau \circ{\rm id} )
\circ (\Delta \otimes \Delta),
$$
where $\tau(v\otimes v')= v' \otimes v$ is the ``flip'' operator on 
$V \otimes V$. This compatibility means indeed that $\Delta$ (or $m$) is
a morphism of algebras (or coalgebras) when a suitable definition of
algebra (or coalgebra) is introduced on $V\otimes V$. 

The application $m$ satisfies some properties that have natural
analogues for $\Delta$, which are  systematically  prefixed with ``co".
For example, the product $m$ is associative, i.e.,   
$$
m\circ (m \otimes {\rm id}) =    m \circ ({\rm id} \otimes m), 
$$
while the coproduct is said to be coassociative
$$
(\Delta \otimes {\rm id}) \circ \Delta = 
({\rm id} \otimes \Delta) \circ \Delta. 
$$    
The product is required to have a unit, and correspondingly the 
coproduct must have  a counit. Algebraically this means that we have 
two linear mappings 
$$
\eta: \extr{C} \longrightarrow V, \qquad 
\epsilon: V  \longrightarrow \extr{C}, 
$$
satisfying
$$
 m \circ ( \eta \otimes {\rm id}) = {\rm id} = m \circ ( {\rm id} 
\otimes \eta),  \qquad  (\epsilon \otimes {\rm id}) \circ \Delta = 
{\rm id} =  ({\rm id} \otimes \epsilon ) \circ \Delta. 
$$

The algebraic structure we have described so far is known as a
bialgebra,  which can be seen as combination of two triplets 
$(V, m, \eta)$ and $(V, \Delta, \epsilon)$ called  algebra and
coalgebra, respectively. 

Hopf algebras are bialgebras characterized besides by the
existence of a linear antimorphism  $\gamma:V \longrightarrow V$ 
verifying 
$$
m \circ (\gamma \otimes {\rm id}) \circ \Delta = \eta \circ \epsilon=
m \circ ({\rm id} \otimes \gamma) \circ \Delta. 
$$
The mapping $\gamma$ is called antipode, and it is easy to show that if
it exists then is unique. 

As examples of this kind of structure we can mention: (finite) group
algebras, the algebra of functions on a (finite, Lie) group,
and enveloping algebras of Lie algebras. All these examples have the
property of being commutative or  cocommutative, i.e., 
$$
m \circ \tau= m, \qquad {\rm or}  \qquad 
\tau \circ \Delta =\Delta,
$$
we can also say that these Hopf algebras are non-deformed 
or ``classical".

\subsect{Quantum algebras and quantum groups}

Quantum groups and quantum algebras are examples of Hopf algebras which
are  neither commutative nor cocommutative. There is a usual
definition of quantum algebra in the sense of Drinfel'd \cite{drinfeld}
and Jimbo \cite{jimbo}, whereas there are several approaches for quantum
groups \cite{faddeev,manin,woronowicz}.

Let ${\fra g}$ be a Lie algebra and ${\cal U}({\fra g})$ its
universal enveloping algebra, which is a ``classical" Hopf algebra 
with coproduct, counit and antipode defined by
$$
\begin{array}{cccc}
&\Delta (X)=1\otimes X+ X\otimes 1, \qquad 
&\Delta (1)= 1\otimes 1;\\[0.2cm]
&\epsilon (X)= 0,\quad \epsilon (1)= 1; \qquad
&\gamma (X)= -X, 
\end{array}
$$
where $X \in {\fra g}$. 

A quantization or deformation of ${\cal U}({\fra g})$ is obtained by
means of a deformed Hopf structure on ${\cal U}_z({\fra g})\equiv 
{\cal U}({\fra g})\hat
\otimes \C[[z]]$, which is the associative algebra of formal power
series in $z$ and coefficients in ${\cal U}({\fra g})$, such that 
$$
{\cal U}_z({\fra g})/z{\cal U}_z({\fra g})\simeq {\cal U}({\fra g})
$$ 
as Hopf algebras (in other words, ${\cal U}_z({\fra g})\to 
{\cal U}({\fra g})$ when $z\to 0$). 

On the other hand, let $G$ be a finite dimensional Lie
group and  ${\fra g}$ its Lie algebra. Let us consider the commutative
and associative algebra of smooth functions of $G$ on $\C$, $Fun(G)$, 
with the usual product of functions (i.e., $(fg)(x)=f(x)g(x), \ f,g
\in  Fun (G),\ x,y \in G$). This algebra has a Hopf structure as
follows 
$$
(\Delta (f)) (x,y)= f (xy),\quad 
\epsilon (f)= f(e), \quad
(\gamma(f))(x)= f(x^{-1}), 
$$
where $f \in Fun (G),\ x,y \in G$, and $e$ is the unit element of $G$.
Note that in general $Fun (G)\otimes Fun (G) \subseteq  Fun (G\times G)$.
When the group is finite the equality is strict, but if $G$ is not
a finite group  $\Delta (f)$ may not belong to 
$Fun (G)\otimes Fun(G)$. This problem can be solved by an
adequate restriction of the space $Fun (G)$.

Incidentally, $Fun (G)$ is the Hopf algebra dual of  ${\cal U}({\fra
g})$  by means of a suitable duality (for more details see, for 
instance, Ref. \cite{charipressley}).

After deformation the above commutative Hopf algebra becomes
non-commutative. On the other hand, $Fun (G)$ is cocommutative if
and only if $G$ is abelian.  

Examples of quantum algebras and quantum groups appear in Section 3. 

\subsect{Bicrossproduct Hopf algebras}

As we mentioned before, a bicrossproduct Hopf algebra can be seen as 
a generalization of the semidirect product of groups
\cite{majid,molnar}. In the following we present the essentials about
this  concept.

Let us start recalling the definition of
$R$--module.  Let $R$ be a unital ring, and ${\cal X}$ a set equipped
with an internal composition law denoted by $+$, and an external
composition law, $\triangleright$, with domain of operators in $R$. We
say that $({\cal X},+,\triangleright)$ is a left
$R$--module if
 
\ \ i)  $({\cal X}, +)$ is an abelian group,  
   
\ ii) the external law ($R\times {\cal X} \to {\cal X}$) satisfies 
\begin{equation}\label{internallow}
\begin{array}{c}
\alpha \triangleright (\beta \triangleright x)= (\alpha \beta) 
\triangleright x, \qquad 
\forall \alpha, \beta \in R,\  \forall x \in {\cal X}, \\[0.2cm]
1 \triangleright x = x,\  \quad\quad \forall x\in {\cal X},
\end{array}
\end{equation}
    
iii) the internal and the external law are compatible in the sense that
\begin{equation} \label{externallow}     
\begin{array}{c}
\alpha\triangleright (x+y)= (\alpha \triangleright x) + 
(\alpha \triangleright y), \qquad 
\forall \alpha \in R,\ \forall x,y \in {\cal X}, \\[0.2cm]
(\alpha + \beta) \triangleright x= (\alpha \triangleright x) + 
(\beta \triangleright x), \qquad
\forall \alpha, \beta \in R,\ \forall x \in {\cal X}. 
\end{array}  
\end{equation}

In the cases of interest for us we will consider the ring associated
with the Hopf algebra $H$, and the set ${\cal X}$ is {\em ab initio} a 
$\extr{C}$--vector space, denoted $V$, hence an abelian group. So, 
we can rewrite the definition of module as follows.

The pair $(V, \rho)$, where $\rho: H \otimes V
\longrightarrow V$ is a linear map, is said to be an $H$--module if 
the external composition law (action) defined by
\be\label{modulo}
h \triangleright v= \rho(h\otimes v)
\ee  
satisfies axioms (\ref{internallow}).
Note that the compatibility conditions (\ref{externallow}) are now
encoded in the linearity  of $\rho$. This mapping is also called a
representation of $H$ on $V$ since it allows to represent the elements
of $H$ by endomorphisms of $V$.

On the other hand, comodules are the dual objects to modules.
The pair $(V, \beta)$, where 
$\beta:  V \longrightarrow H \otimes V \  (\beta(v)=v^{(1)} 
\otimes v^{(2)})$ is
a linear map, is said to be an $H$--comodule if the 
``external decomposition'' law (coaction) defined by $\beta$ satisfies
\be\label{comodulo}
\begin{array}{c}
(v^{(1)} \leco ) \otimes v^{(2)}= 
\Delta(v^{(1)})\otimes v^{(2)}, \qquad \forall v \in V,\\[0.2cm]
\epsilon(v^{(1)}) \otimes v^{(2)}= v, \qquad\qquad\qquad\qquad
\forall v \in V, 
\end{array}
\ee
where we have written $\beta(v)= v \leco$
to make the notation more symmetric between actions and coactions.
Similarly to $\rho$ the mapping $\beta$ is called corepresentation.

Remark that for the last two definitions and do not take into account
the whole Hopf algebra structure. Thus, only the algebra (coalgebra)
sector is used for modules (comodules).   

When a Hopf algebra $H$ acts on an algebra $A$ it is
natural to demand some compatibility of the action with the algebraic
structure. So, we say that $A$  is a right $H$--module algebra if it is
an $H$--module and the action satisfies 
$$
 (a a')\triangleleft h= (a\triangleleft h_{(1)})
(a'\triangleleft h_{(2)}), 
\qquad  1 \triangleleft h= \epsilon(h),
$$
where $\Delta(h)=\sum  h_{(1)}\otimes h_{(2)}$.

A similar situation happens when $H$ coacts on a coalgebra $C$. Then, 
it is said that $C$ is a left $H$--comodule coalgebra if it is an
$H$--comodule and the coaction ($\leco$) satisfies
$$
\begin{array}{c}
  c^{(1)} \epsilon(c^{(2)})= 1_H \epsilon(c), \\[0.2cm]
c^{(1)} \otimes {c^{(2)}}_{(1)} \otimes {c^{(2)}}_{(2)} =
{c_{(1)}}^{(1)} {c_{(2)}}^{(1)} \otimes {c_{(1)}}^{(2)} 
\otimes {c_{(2)}}^{(2)},
\end{array}
$$
where $c\leco=c^{(1)}\otimes c^{(2)}$.
 
It is also possible to define module coalgebras and  
comodule algebras (both with right and left-handed versions, of course) 
but we shall not need it.

The introduction of module algebras allows to go a step beyond the 
direct sum of algebras. If $A$ is a right $H$--module algebra 
we can define an algebra structure on the tensor product $H \otimes A$ 
by means of the composition law 
$$
(h\otimes a)(h'\otimes a')= h h'_{(1)} \otimes 
(a \triangleleft h'_{(2)}) a'.
$$
It is immediate to check that the new algebra has $1_H \otimes 1_A$ as 
unit element. This structure is called the semidirect product $H \RIMO
A$.

By duality, if we start with the left $H$--comodule coalgebra 
$C$ we can construct a coalgebra structure on the tensor product $C
\otimes H$ defining the coproduct as
$$
\Delta(c\otimes h)= c_{(1)} \otimes {c_{(2)}}^{(1)} h_{(1)} 
\otimes   {c_{(2)}}^{(2)} \otimes h_{(2)},
$$
and taking $\epsilon_C \otimes \epsilon_H$ as antipode. The resulting
coalgebra is the semidirect product $C \LECO H $.

Now, let us change the notation to consider simultaneously 
two Hopf algebras $K$ and $L$, such that $L$ is a right $K$--module 
algebra and $K$ is a left $L$--comodule  coalgebra. According to the 
above semidirect product constructions  $K\otimes L$ is equipped with 
a structure of algebra ($K \RIMO L$) and other of coalgebra  ($K \LECO
L$).  The following five compatibility conditions \cite{majid} 
$$
\begin{array}{c}
\epsilon(l \triangleleft k)= \epsilon(l) \epsilon(k),  \quad
\Delta(l\triangleleft k)= (l_{(1)} \triangleleft k_{(1)})
{k_{(2)}}^{(1)}  \otimes l_{(2)} \triangleleft {k_{(2)}}^{(2)},
\\[0.2cm] 
1 \leco = 1 \otimes 1, \quad  (kk') \leco= 
(k^{(1)}\triangleleft k'_{(1)}){k'_{(2)}}^{(1)} \otimes
k^{(2)}{k_{(2)}'}^{(2)},\\[0.2cm]
{k_{(1)}}^{(1)} (l \triangleleft k_{(2)}) \otimes
{k_{(1)}}^{(2)}=  (l \triangleleft k_{(1)}) {k_{(2)}}^{(1)} \otimes
{k_{(2)}}^{(2)} ,
\end{array}
$$
are sufficient conditions to guarantee that both structures fit 
adequately to form a bialgebra with antipode: the right-left
bicrossproduct Hopf algebra  $K \bicross L$. The antipode is given by
$$
\gamma (k\otimes l)= (1\otimes \gamma (k^{(1)}l) )(\gamma (k^{(2)})
\otimes1).  
$$
In analogy with the classical case we shall refer to $L$ as the  
kernel of the bicrossproduct.  
For our purposes we are interested in the case in which the kernel is
commutative.

Since $K$ and $L$ generate $K\bicross L$ we can 
construct a basis of $K\bicross L$ using bases of $K$ and $L$.

When $K$ and $L$ are the universal enveloping algebras of  Lie algebras 
$\fra k$ and $\fra l$, respectively, the right action of
$K$ on $L$  is  given by means of the Lie commutators, i.e.,  
$l\triangleleft {k}= [l,{k}], \ l\in L,\ k\in {\fra k}$.

The left-right version is constructed in a similar way. 
In this case one considers the right $\cal K$--comodule coalgebra 
$\cal L$ and the left $\cal L$--module algebra $\cal K$. The new product
and coproduct on $\cal K \otimes L$ are defined by
$$
\begin{array}{c}
(\kappa \otimes \lambda)(\kappa' \otimes \lambda')= 
\kappa(\lambda_{(1)} \triangleright \kappa') \otimes \lambda_{(2)}
\lambda', \\[0.2cm]
\Delta(\kappa \otimes \lambda)= (\kappa_{(1)}  \otimes
{\lambda_{(1)}}^{(1)}) \otimes 
(\kappa_{(2)} {\lambda_{(1)}}^{(2)} \otimes \lambda_{(2)}).
\end{array}
$$
The unit and counit are as in the right-left case. The compatibility 
conditions read off as 
$$
\begin{array}{c}
\v(\l\act \k)=\v(\l)\v(\k),\\[0.2cm]
\Delta(\l\act \k)\equiv (\l\act \k)_{(1)}\otimes (\l\act \k)_{(2)}
=({\l_{(1)}}^{(1)}\act \k_{(1)})\otimes {\l_{(1)}}^{(2)}(\l_{(2)}
\act \k_{(2)}),\\[0.2cm]
\reco (1)\equiv 1^{(1)}\otimes 1^{(2)}= 1\otimes 1,\\[0.2cm]
\reco (\k\k')\equiv (\k\k')^{(1)}\otimes (\k\k')^{(2)}
={\k_{(1)}}^{(1)}\k'^{(1)}\otimes {k_{(1)}}^{(2)}(\k_{(2)}
\act \k'^{(2)}) ,\\[0.2cm]
{\l_{(2)}}^{(1)}\otimes(\l_{(1)}\act \k){\l_{(2)}}^{(2)}=
{\l_{(1)}}^{(1)}\otimes {\l_{(1)}}^{(2)}(\l_{(2)}\act \k).
\end{array}
$$
The left-right bicrossproduct structure is denoted by $\cal K \LR L$.

In the finite dimensional case it is easy to show that
$(K\bicross L)^*= K^* \LR L^*$. For the cases we are interested in, 
although they are infinite dimensional, a similar result holds,
provided  that ``duality'' is changed for ``dually paired algebras''
(see Ref. \cite{majid} or \cite{charipressley}). 

On the other hand, one can prove that given a  bicrossproduct Hopf 
algebra, $H=K \LR L$, with dual $H^*=K^* \RL L^*$, a nondegenerate dual
pairing between $H$ and
$H^*$ can be defined in terms of nondegenerate pairings $\langle \cdot, 
\cdot \rangle_1$ and $\langle \cdot, \cdot \rangle_2$ for the pairs
$(K,K^*)$ and $(L, L^*)$, respectively, by
$$
\langle kl, \kappa \lambda \rangle= \langle k, \kappa 
\rangle_{1}  \langle l, \lambda \rangle_{2}. 
$$
An immediate consequence of this statement, that we shall  use later, 
is  that with the pairings defined above, if  $\{ k_m\}$ and
$\{\kappa^m\}$ are dual basis for $K$ and $K^*$, and
$\{l_n\}$ and $\{\lambda^n\}$ are dual basis for $L$ and $L^*$,
respectively, then $\{k_m
l_n\}$ and
$\{\kappa^m \lambda^n\}$ are dual basis for $H$ and $H^*$.

\sect{Induced representations}

In the theory of representations of Hopf algebras the symmetry
played by the algebra and coalgebra structures is broken.
On the one hand, the algebra structure of a Hopf algebra, $H$, leads to 
a ring structure on $H$ and, hence, to (in general non-commutative)
module theory, but on the other hand, the coalgebra structure allows a
tensor product of $H$--modules turning this category into a monoidal
one.

Induced representations are precisely extensions of scalars from
the point of view of module theory. Effectively, let us consider the
unital  associative algebra $A$ as a ring with unit, let $B$ be a
subalgebra of $A$ containing the unit and $V$ a right $B$--module. The
algebra $A$ can be considered as a left $B$--module (by means of left
regular translations), and therefore  the tensor product (on $B$)  $V
\otimes_B A$ makes sense. In this last $B$--module, $V \otimes_B A$, we
can extend the scalars to $A$, and then we say that $V^\uparrow=V\otimes_B
A$ is the $A$--module induced from the $B$--module $V$. 

A similar construction can  be carried out by replacing  $A$ by
its linear dual  $A^*$, and looking at it as the right $B$--module
associated with the left regular action on $A$. In the literature this
construction is referred to as coinduced  representations \cite{dixmier}
or produced representations \cite{higman}.
 
It is worthy to note that the terminology about representations in 
quantum group literature is a bit confusing. So, terms like
``representation" and ``induced representation'' have been also used to
denote  ``corepresentation'' and  ``induced corepresentation'' 
\cite{ibort,giller,maslanka}. Hoping not introduce
more confusion  we will speak of induced representations in the sense
of the preceding paragraph.

The ideas introduced above allow us to develop 
the construction of the induced representations for two interesting
physical examples, corresponding to two non-equivalent deformations
of the one-dimensional Galilei algebra. 

\subsect{Standard quantum $(1+1)$ Galilei algebra}

The standard quantum $(1+1)$ Galilei algebra for which we calculate
the induced representations is a contraction of the $\kappa$-Poincar\'e
in $(1+1)$ dimensions \cite{azcarraga1}.

\subsubsect{Algebraic structure}
The Hopf algebra structure of the standard quantum $(1+1)$ Galilei
algebra, $U_\omega [{\fra g}(1,1)]$, is defined by  
\begin{eqnarray*}
&[H, K]= - P , \quad [P, K]= \omega P^2 , \quad [H, P]=0 ;  &
\nonumber \\ &\Delta H = H \otimes 1 + 1 \otimes H , \quad 
\Delta X = X \otimes 1 + \exp (- 2 \omega H) \otimes X , 
\ \   X \in \{P, K \} ;   &  \nonumber \\
&\epsilon(X)= 0  , \quad X \in \{H, P, K\};  \\
& \gamma(H)= -H,  \quad  \quad  \gamma(X)= - e^{2 \omega H} X,  
\quad \quad X \in \{P, K\}. & \nonumber  
\end{eqnarray*} 

The Hopf algebra  $U_\omega [{\fra g}(1,1)]$ has a  bicrossproduct 
structure \cite{azcarraga1} given by $U_\omega [{\fra g}(1,1)]=
U [\R] \bicross U_\omega [{\fra t}_2]$, where 
$U [\R]= \langle K \rangle$  (``boost sector'') and 
$U_\omega [{\fra t}_2]$ is the deformed subalgebra  generated by $P$ and 
$H$ (``translation sector''). The right action of
$U [\R]$ on $U_\omega [{\fra t}_2]$ is given by  
$$
P \triangleleft K= \omega P^2, \qquad H \triangleleft K = -P,
$$
and the left coaction of $U_\omega [{\fra t}_2]$ on $U [{\R}]$ is
$$
 K \leco= e^{-2 \omega H} \otimes K.
$$

The dual algebra $F_\omega [G(1,1)]$ is generated by the local
coordinates $v, x, t$. The commutators, coproduct, counit and 
antipode are given by \cite{azcarraga1}
\begin{eqnarray*}
&[t, x]= -2 \omega x , \quad [x,v] = \omega v^2  , 
\quad [t, v]= - 2 \omega v ;  & \nonumber \\[0,2cm] 
&\Delta t = t \otimes 1 + 1 \otimes t , \quad  
\Delta x = x \otimes 1 + 1 \otimes x - t \otimes v , \quad  
\Delta v = v \otimes 1 + 1 \otimes v  ; &  \nonumber \\[0,2cm] 
& \epsilon(f)= 0,\qquad f\in \{v, t, x\} ; &  \\[0,2cm] 
& \gamma (v)= -v,  \qquad  \gamma (x)= -x - t v,  \qquad   
\gamma (t)= -t. & \nonumber 
\end{eqnarray*}
The bicrossproduct structure  $F_\omega [G(1,1)]= \langle v \rangle
\LR \langle  x, t \rangle$ is encoded in the left action
$$
 x \triangleright v= \omega v^2, \qquad t \triangleright v= -2 \omega v,
$$
and right coaction 
$$
\reco x= x \otimes 1 - t \otimes v, \qquad \reco t= t \otimes 1.
$$

The duality pairing between the Hopf algebras 
$U_\omega [{\fra g}(1,1)]$ and $F_\omega [G(1,1)]$
is given explicitly by
\begin{equation}
\langle K^m P^n H^p, v^q x^r t^s \rangle = 
m! n! p! \,  \delta^m_q \delta^n_r \delta^p_s.\label{pairing1} 
\end{equation}

\subsubsect{Induced representations} 

The representations of $U_\omega [{\fra g}(1,1)]$
are induced using (left--)characters of the
translation sector  
\be \label{leftcharacteres}
1 \dashv P^n H^p = (i a)^n (i b)^p,\qquad  n,p \in \extr{N},
\ee
where $\dashv$ stands for the action of $U_\omega [{\fra t}_2]$ on 
$\extr{C}$. These induced representations have as carrier space
$\extr{C}^\uparrow$ the space
${\rm Hom}_{U_\omega [{\fra t}_2]}(U_\omega[ {\fra g}(1,1)], 
\extr{C})$, which  is contained in the space  
${\rm Hom}_{\extr{C}}(U_\omega[ {\fra g}(1,1)], \extr{C})= 
U_\omega[{\fra g}(1,1)]^*$. This later can be identified with
$F_\omega[G(1,1)]$. Since the elements of $\extr{C}^\uparrow$ are 
$U_\omega [{\fra t}_2]$--morphisms they are characterized by the
``equivariance condition''  
\begin{equation}\label{equivariancecondition} f(X P^nH^p)= f(X) \dashv
P^n H^p. \end{equation}
A generic element of $F_\omega[G(1,1)]$ will  be
\begin{equation} 
f= f_{q,r,s} v^q x^r t^s, \label{funcion}
\end{equation}
then, using the pairing (\ref{pairing1}) and imposing the
equivariance condition (\ref{equivariancecondition}), in  order to have
$f$ contained in $\extr{C}^\uparrow$, we get 
\be \label{funcion1}
 q! r! s! f_{q,r,s} = \langle f, K^q P^r H^s \rangle = 
\langle f, K^q \rangle \dashv P^r H^s = q! f_{q,0,0} (ia)^r (ib)^s.
\ee
Introducing this last relation (\ref{funcion1}) in expression
(\ref{funcion}) we obtain that  $\extr{C}^\uparrow$ is the subspace of
$F_\omega [G(1,1)]$  whose elements are of the form $\phi(v)  e^{i a x}
e^{i b t}$. 

Let us consider the basis $v^m e^{iax}e^{ibt}$ of ${\C}^\uparrow$. The
action of the elements $X$ of $U_\omega [{\fra g}(1,1)]$ on it  will be
given in terms of the ${\C}$--numbers $[X]^m_{q,r,s}$ by means of the
expression 
\begin{equation} 
 v^m e^{iax} e^{ibt} \dashv X = [X]^m_{q,r,s} v^q x^r t^s,
\label{cnumeros} 
\end{equation}
where $\dashv$ also denotes the action of $U_\omega [{\fra g}(1,1)]$
on $\extr{C}^\uparrow$. The evaluation of $[X]^m_{q,r,s}$ is made using 
the pairing (\ref{pairing1}) 
$$
q! r! s! [X]^m_{q,r,s} = 
\langle v^m e^{iax} e^{ibt} \dashv X , K^q P^r H^s \rangle = 
\langle v^m e^{iax} e^{ibt},  X K^q P^r H^s \rangle. 
$$
So, the computation of  $[X]^m_{q,r,s}$ has been reduced to the 
problem of writing the monomial $X K^q P^r H^s$ in the ``normal
ordering'' defined by the above basis of ${\C}^\uparrow$.
When  $X=K$ the task is trivial, for the other two generators, $P$
and $H$, we use the following results whose proof is made by induction
$$
\begin{array}{ll}
P\LLL K^k = k! \omega^k  P^{k+1},& 
 P K^q= {\displaystyle \sum_{k=0}^q} \frac{q!}{(q-k)!} \omega^k 
 K^{q-k} P^{k+1}, \\  
H\LLL K^{k+1} =- k! \omega^k P^{k+1}, & 
H K^q= K^q H - {\displaystyle \sum_{k=0}^{q-1}} \frac{q!}{(k+1) 
(q-k-1)!} \omega^k K^{q-k-1} P^{k+1}. \\
\end{array}
$$
Thus, we get that
\begin{eqnarray}  \label{caccion}
 &  q! r! s! [K]^m_{q,r,s}  = m! \delta^m_{q+1} (ia)^r (ib)^s,  & 
\nonumber \\[0.2cm]
 &  q! r! s! [P]^m_{q,r,s}= m! {\displaystyle \sum_{k=0}^q }
  \frac{q!}{(q-k)!}  \omega^k \delta^m_{q-k} (ia)^{k+1+r} (ib)^s,  
  & \nonumber \\[0.2cm]
 &  q! r! s! [H]^m_{q,r,s}= m! \delta^m_q (ia)^r (ib)^{s+1}\qquad\qquad
\qquad\qquad&\\[0.2cm] 
&\qquad\qquad\qquad\qquad\qquad - m! {\displaystyle
\sum_{k=0}^{q-1}} \frac{q!}{(k+1) (q-k-1)!}\delta^m_{q-k-1} \omega^k
(ia)^{k+1+r} (ib)^s.  & \nonumber 
\end{eqnarray}
Now substituting expressions (\ref{caccion}) in (\ref{cnumeros}) 
we obtain the desired action of the generators on the basis of 
$\extr{C}^\uparrow$. Finally, in order to have meaningful
expressions for the actions of the generators of 
$U_\omega [{\fra g}(1,1)]$ it is necessary ``to complete" the space
$\extr{C}^\uparrow$. In consequence, we shall work with the space of
formal series in $v$, $\extr{C}[[v]]$.

Summarizing, the induced representations of $U_\omega [{\fra g}(1,1)]$
determined by the (left--)characters of the 
translation sector (\ref{leftcharacteres}) have as support space the 
space  $\extr{C}[[v]]$. The explicit form of
these representations is  
\begin{equation}\label{representation1}
{\begin{array}{l}
\phi(v)\dashv K = \phi'(v),  \\[0.2cm]
\phi(v)\dashv P=  \phi(v) \, \frac{ia}{1 - \omega ia v},  \\[0.2cm]
\phi(v)\dashv H=  \phi(v) \,  [ib + \frac{1}{\omega} 
\ln (1-i \omega a v)]. 
\end{array}}
\end{equation}
The representation is labeled by two real parameters $a$ and $b$, 
however, by the transformation $H\to H-ib$ the coefficient $b$ vanishes.
Therefore, the representations labeled by $(a,b)$ are pseudoequivalent 
(i.e., equivalent up to a phase) to those with  $(a,0)$. 

It is worthy to note that in the limit $\omega \to 0$ we recover a
unitary irreducible representation of the nondeformed Galilei group
provided that $a$ is a real parameter. Effectively, taking the
limit $\omega \to 0$ in (\ref{representation1}) we get the infinitesimal
action of the infinitesimal generators of the group 
\begin{equation}\label{representation2}  {\begin{array}{ll}
\phi(v)\dashv K = \phi'(v), \qquad & K=d/dv, \\[0.2cm] 
\phi(v)\dashv P=  ia \phi(v) ,  \qquad & P=ia,\\[0.2cm] 
\phi(v)\dashv H= -iav  \phi(v), \qquad & H=-iav .  
\end{array}}
\end{equation}
On the other hand, the unitary irreducible representation up to a phase 
obtained by Mackey's method that correspond to the above one ($m=0$,
$C\neq 0$) is \cite{miguelangel}
\be 
\label{representation3} 
{\cal U}_p(t,x,v)\psi)(\xi)=e^{i(px-\xi t)}\psi(\xi +vp),\qquad p\in
\R^*.
\ee
If one computes the  infinitesimal action associated with the above
representation (\ref{representation3}) it coincides with 
(\ref{representation2}) when $p=a$ and after the variable change 
$v\leftrightarrow \xi=va$. Obviously, the functions to be consider in 
the limit will be of integrable square. 

So, the space of formal series in $v$ is too large for the study of the
unitarity as we have just mentioned. For that, we can reduce the
support space by considering the space of polynomials in
$v$. Let us define 
$$
v_k=(e^{iax})^{-k} v(e^{iax})^k= \frac{v}{1-ki\omega av},\qquad k\in
\extr Z . 
$$ 
The action (\ref{representation1}) in the space of the
polynomials in $v_k$ is as follows
\begin{eqnarray}\label{representation4} 
(v_k)^n\dashv K& =& n(v_k)^{n-1}{(1+ki\omega av_k)^2}, 
\nonumber \\[0.2cm] 
(v_k)^n\dashv P &= &  ia(v_k)^{n}{(1+ki\omega av_1)},  
\\[0.2cm] 
(v_k)^n\dashv e^{\omega H} &= &
ia(v_k)^{n}{(1-ki\omega av_0)},\nonumber  
\end{eqnarray}
where $v_0=v$. We see that the representation (\ref{representation4}) 
is reducible but non completely reducible, and 
$\C \oplus {\extr P} [v_0] \oplus {\extr P} [v_1]$ determines an
irreducible subspace for this representation, where ${\extr P} [v_k]=v_k
\C [v_k]$ and $\C [v_k]$ is the space of polynomials in $v_k$. An open
problem is to construct a Haar measure in such a way that this
representation becomes unitary, and in the limit we can recover  the
space of integrable square functions. A solution of this problem for the
quantum  $(1+1)$ extended Galilei group is given in Ref.
\cite{italianos}.

\subsect{Non-standard quantum $(1+1)$ Galilei algebra}
 
This non-standard quantum $(1+1)$ Galilei algebra is a contraction 
\cite{azcarraga1} of the non-standard Poincar\'e algebra
\cite{ballesteros1,ballesteros2}.

\subsubsect{Algebraic structure} 

The structure of $U_\rho [{\fra g}(1,1)]$ is given by 
\begin{eqnarray*}
&[H, K]= - \frac{1}{4\rho} (1- \exp(-4 \rho P)) , 
\quad [P, K]= 0 , \quad [H, P]=0  ;  & \nonumber \\[0.2cm]
&\Delta P = P \otimes 1 + 1 \otimes P , \quad 
\Delta X = X \otimes 1 + \exp (- 2 \rho P) \otimes X , 
\ \  X \in \{H, K \}  ;  & \nonumber \\[0.2cm]
&\epsilon(X)= 0  , \quad X \in \{H, P, K\};  & \\[0.2cm] 
&\gamma (P)= -P,  \quad \quad  \gamma (X)= - e^{2 \rho P} X,  
\quad \quad X \in \{H, K\}. &  
\end{eqnarray*} 

The quantum group $F_\rho [G(1,1)]$ is determined by 
\begin{eqnarray*}
& [t, v]= 0  , \quad [x, v]= - 2 \rho v , 
\quad [t,x] = 2 \rho t ; & \nonumber \\[0.2cm]
&\Delta t = t \otimes 1 + 1 \otimes t , \quad 
\Delta x = x \otimes 1 + 1 \otimes x - t \otimes v , \quad
\Delta v = v \otimes 1 + 1 \otimes v  ; & \nonumber \\[0.2cm]
& \epsilon(f)=0 , \quad f \in \{t, x, v\}; &\\[0.2cm]
&\gamma (v)= -v,  \quad \quad  \gamma (x)= -x - t v,  
\quad \quad  \gamma (t)= -t.& \nonumber
\end{eqnarray*}
We can make similar considerations to those made for the standard case
about the  bicrossproduct structure, except that now  the duality 
between both algebras is given by  
\begin{equation} \label{pairing2}
\langle
K^m H^n P^p, v^q t^r x^s \rangle = 
m! n! p! \delta^m_q \delta^n_r \delta^p_s .
\end{equation}
Note also that now the order of $H,P$ and $t,x$ has been changed with
respect to the order taken in (\ref{pairing1}).

\subsubsect{Induced representations} 

Let us consider the following representation of  the translation
sector $L$ on $\extr{C}$  
\begin{equation}
1 \dashv H^n P^p = (i b)^n (i a)^p, \quad n,p \in \extr{N}.
\end{equation}
The support space of the induced representation, denoted by
$\extr{C}^\uparrow$, is the subspace of  $F_\rho [G(1,1)]$  whose
elements are like
$\phi(v)  e^{i bt} e^{i ax}$.  This  subspace is isomorphic to
$\extr{C}[[v]]$. The explicit action of the generators of $U_\rho
[{\fra g}(1,1)]$ over the elements of $\extr{C}[[v]]$ is
\begin{equation} {\begin{array}{l}
\phi (v) \dashv K =  \phi'(v),   \\[0.2cm]
\phi (v) \dashv P=   \phi (v) \, i a,   \\[0.2cm]
\phi (v) \dashv H=   \phi (v) \, [i b + 
\frac{1}{4\rho}(1- e^{-4ia \rho}) v].  
\end{array} }\label{representation5}
\end{equation}
The computation of this representation is based on the following
result
$$
\begin{array}{l}
H\LLL K^k = H \delta ^k_0 - \frac{1}{4 \rho} (1 - e^{- 4 \rho P}) 
\delta^k_1, \\[0.2cm]
  H K^q= K^q H - \frac{1}{4 \rho} q  K^{q-1} (1 - e^{- 4 \rho P}).
\end{array}
$$

Similarly to the above case the representation labeled by $(a,b)$ is
equivalent to that with $(a,0)$. 

The irreducibility of the representation follows from the fact that
$K$ and $H$ can be interpreted as ladder operators acting on the space
of polynomials in $v$. The above result is not all surprising if one
takes into  account that the algebra (only the algebra, not the whole
Hopf structure) $U_\rho [{\fra g}(1,1)]$
contains the oscillator algebra.

\sect{Conclusions}

We have constructed induced representations of two quantum groups,
and seen that only the algebra structure has been relevant in
our procedure. The coalgebra structure helps in the computation of
some expressions but is not essential. However, the coalgebra structure
is crucial to allow the tensor product of representations of the
Hopf algebra.

Both examples presented here have a  bicrossproduct structure,
which provides technical facilities,  for example, in evaluating
pairings, nevertheless it has not been essential in the induction
process.

The mechanism of induced representations looks to be a systematic 
way of construction of representations, while some times in the
literature the construction of representations or corepresentations has
been made by {\em ad hoc} procedures. 

There are some open problems to establish a complete theory of induced
representations of quantum groups. We can mentioned, for instance, the
definition of equivalence criteria of representations; the
irreducibility of the representations, that is, to know the conditions
to construct irreducible representations; the unitarity of the induced
representations, and if this procedure allows to obtain all the
irreducible representations. A solution for the unitarity problem is
connected with the construction of a quantum analogue of the Haar
measure. Work on these questions is in progress, and the results will
be published elsewhere.


\section*{Acknowledgments}  This work has been partially supported by
DGES of the  Ministerio de Educaci\'on y Cultura of Spain under Project
PB95-0719,  and the Junta de Castilla y Le\'on (Spain).


\end{document}